\documentclass[12pt]{article}
\usepackage[dvips]{graphics}
\usepackage{amsthm,amsfonts,amsmath, amssymb}
\usepackage{alltt,makeidx,enumerate}
\usepackage{epsfig}
\usepackage{subfig}
\usepackage{color}
\usepackage{soul}
\allowdisplaybreaks

\newtheorem{lem}{Lemma}[section]
\newtheorem{thm}[lem]{Theorem}

\newtheorem{example}[lem]{Example}

\newenvironment{prueba}{\noindent\textit{Proof}:}{  \qed\\\indent}

\newcommand{\dif}{\mathrm{d}}

\newcommand{\limite}[2]{\lim_{#1\rightarrow #2}}
\newcommand{\deru}[1]{#1^\prime}

\newcommand{\norm}[1]{\left\Vert#1\right\Vert}
\newcommand{\abs}[1]{\left\vert#1\right\vert}
\newcommand{\set}[1]{\left\{#1\right\}}
\newcommand{\mc}[1]{\mathcal{#1}}
\newcommand{\Real}{\mathbb R}
\newcommand{\eps}{\varepsilon}
\newcommand{\To}{\rightarrow}

\newcommand{\B}{\mathcal{B}}
\newcommand{\ts}[1]{\mathbf{#1}}

\newcommand{\pd}[2]{\frac{\partial#1}{\partial#2}}
\newcommand{\td}[2]{\frac{\dif#1}{\dif#2}}
\newcommand{\tdd}[2]{\frac{\dif^2#1}{\dif#2^2}}

\newcommand{\be}{\begin{equation}}
\newcommand{\ee}{\end{equation}}

\begin{document}

\title{Cavitation of a spherical body under mechanical and self gravitational 
forces}

\author{Pablo V. Negr\'on--Marrero\thanks{pablo.negron1@upr.edu}\\
Department of Mathematics\\
University of Puerto Rico\\
Humacao, PR 00791-4300\\\and Jeyabal
Sivaloganathan\thanks{masjs@bath.ac.uk}\\
Department of Mathematical Sciences\\University of Bath, Bath\\BA2
7AY, UK }

\date{}

\maketitle

\begin{abstract}
In this paper we look for minimizers of the energy functional for 
isotropic 
compressible elasticity taking into consideration the effect of a 
gravitational field induced by the body itself. We consider the displacement 
problem in which the outer boundary of the body is subjected to a Dirichlet 
type boundary 
condition.  For a spherically symmetric body 
occupying the unit ball $\B\in\Real^3$,  the minimization is done within the 
class of radially symmetric 
deformations. We give conditions for the existence of such minimizers, for 
satisfaction of the Euler--Lagrange equations, and show that for large 
displacements the minimizer must develop a cavity at the centre. A numerical 
scheme for approximating these minimizers is given together with some 
simulations that show the dependence of 
the cavity radius and minimum energy on the displacement and mass density of 
the body.
\end{abstract}

\section{Introduction}
The study of the shape of self gravitating bodies is extensive and dates back 
to the time of Newton itself. It is well known that depending on the density of 
a dying star, there are several possibilities for the resulting object: white 
dwarf, neutron star, black hole, etc. The case of a black hole forming is also 
referred to as \textit{gravitational collapse}. The literature on these 
phenomena is extensive and we refer to \cite{CaLe2012} and \cite{JiKoChGo2019} 
for a historical account.

In this paper we consider the problem of a self gravitating spherical body. 
Apart 
from its apparent ``simplicity'', this problem plays an important role on the 
study of the more complex phenomena described above. The proposed 
variational model combines both mechanical and gravitational responses, the 
mechanical part based on a model from nonlinear elasticity which allows for the 
characterization of large deformations. Under certain mathematically 
physical conditions, the extrema of the corresponding energy functional, can be 
characterized via the Euler--Lagrange equations. This combined model has been 
used by \cite{BeSc2003}, \cite{CaLe2012} and \cite{JiKoChGo2019} among others. 
In \cite{BeSc2003} the existence of solutions to the Euler--Lagrange equations,
with a zero dead load boundary condition on the outer boundary of the body, is 
established via the implicit function theorem and is valid for ``small bodies'' 
of arbitrary shape. 

The work in \cite{JiKoChGo2019} is for spherically symmetric deformations with a 
zero dead load boundary on the outer boundary condition as well, and combines 
asymptotic analysis 
with numerics to get results for varying densities and reference configuration 
body radius. They used a stored energy function of the form 
\begin{equation}\label{se_collapse}
 W(\ts{F})=\frac{\mu}{2}\left(\norm{\ts{F}}^2-3-2f_\alpha(\det\ts{F})\right) 
+\frac{\beta}{2}\left(\det\ts{F}-1\right)^2,
\end{equation}
where $\alpha\ge0$, $\beta$ and $\mu$ are positive constants, and 
$f_\alpha(d)=\ln(d)-\alpha 
d^{-\alpha}(d-1)^4$. This material corresponds to a ``soft'' compressible 
material for $\alpha=0$ or small, and to a ``strong'' compressible material 
otherwise. For constant reference configuration density $\rho_0$, the authors 
in \cite{JiKoChGo2019} show numerically that for $\alpha$ small there exists a 
critical density $\rho_0^*$ such that if $\rho_0\le\rho_0^*$, then 
the Euler--Lagrange equations (cf. \eqref{eqn2.13}) can have multiple 
solutions, most of them unstable, while if $\rho_0>\rho_0^*$, then there are no 
solutions which could be interpreted as gravitational collapse. Moreover for 
$\alpha$ large, there are 
solutions for all densities $\rho_0$, which appear to be unique.

By adapting the techniques in \cite{Ba77} for polyconvex stored energy 
functions, the authors in \cite{CaLe2012} show the existence of minimizers for 
the resulting energy functional, now for large deformations and arbitrary 
bodies, and for both, zero dead load and  displacement boundary conditions. The 
stored energy functions used in \cite{CaLe2012} could be classified as 
corresponding to  
``strong'' compressible materials (cf. eqns. (13) and (14) in \cite{CaLe2012}).

In this paper we look for minimizers of the energy functional for isotropic 
compressible elasticity and taking into consideration the effect of a 
gravitational field induced by the body itself. We consider the displacement 
problem in which the outer boundary of the body is subjected to a Dirichlet 
type boundary 
condition.  For a spherically symmetric body 
occupying the unit ball $\B\in\Real^3$ and with radially symmetric mass 
density,  the minimization is done within the class of radially symmetric 
deformations. Contrary to previous works, the deformations we consider 
belong to $W^{1,p}(\B)$ with $p<3$, and thus may develop singularities. For the 
particular case of radially symmetric deformations, we study the occurrence or 
initiation of a cavitation at the centre of the ball and its dependence on the 
boundary displacement and gravitational related constants.

In Section \ref{sec:2} we introduce the basic model, with energy functional and 
admissible function space, for radial deformations (cf. \eqref{eqn2.8}) of a 
spherically symmetric body. These deformations are characterized by a function 
$r\,:\,[0,1]\To[0,\infty)$, the Dirichlet boundary condition taking the form 
$r(1)=\lambda$. After this we show in Section \ref{sec:3} that under certain 
growth conditions on the stored energy function (cf. \eqref{hdet_growth} with 
H1--H3) and for any reference configuration density function $\rho_0$ that is 
bounded, nonnegative and bounded away from zero, a minimizer of the energy 
functional \eqref{eqn2.10}--\eqref{potfunc} exists over the admissible set 
\eqref{alambda}. Under the additional constitutive assumption \eqref{derPhibdd}, 
these minimizers satisfy the Euler--Lagrange equations \eqref{eqn2.13} where 
either $r(0)=0$ or $r(0)>0$ (\textit{cavitation}) with zero Cauchy 
stress at the origin (cf. \eqref{cavBC}). In Section \ref{sec:4} we show that 
for $\lambda$ sufficiently 
large, these minimizers must satisfy $r(0)>0$. This result is an adaptation to 
the problem with self gravity of a similar result in \cite{Siva1991}  
for compressible inhomogeneous materials.

In Section \ref{sec:5} we collect several results for $\lambda$ small where the 
minimizers must have the centre intact. In addition we show in Theorem 
\ref{nocav_ls} that any minimizer which leaves the centre intact must have 
strains at the origin less than the critical boundary displacement 
corresponding to an isotropic material made of the material at the centre of  
the original body. Once again this result is an adaptation to 
the problem with self gravity of a similar result in \cite{Siva1991}  
for compressible inhomogeneous materials. 

Finally in Section \ref{sec:6} we 
present a numerical scheme for the computation of the minimizers of our energy 
functional. This method is based on a combination of a gradient flow iteration 
which works as a predictor, together with a shooting method to solve the  
EL-equations, that works as a corrector. For constant reference configuration 
densities we present several simulations that show the dependence of 
the cavity radius and minimum energy on the displacement $\lambda$ and density 
$\rho_0$.

\section{Problem formulation}\label{sec:2}
Consider a body which in its reference configuration occupies
the region
\begin{equation}\label{eqn2.1}
\B=\{ \ts{x}\in \Real ^3\,|\, \norm{\ts{x}}<1\},
\end{equation}
where $\norm{\cdot}$  denotes the Euclidean norm. Let
$\ts{u}:\B\rightarrow\Real^3$ denote a deformation of the body and
let its \textit{deformation gradient} be
\begin{equation}\label{eqn2.2}
\nabla\ts{u}(\ts{x})=\td{\ts{u}}{\ts{x}}(\ts{x}).
\end{equation}
For smooth deformations, the requirement that $\ts{u}(\ts{x})$ is
locally \textit{invertible and preserves orientation} takes the
form
\begin{equation}\label{eqn2.3}
\det \nabla\ts{u}(\ts{x})>0,\quad\ts{x}\in\B.
\end{equation}
Let $W:M_+^{3\times 3}\rightarrow\Real$ be the \textit{stored
energy function} of the material of the body where $M_+^{3\times
3}=\{ \ts{F}\in M^{3\times 3}\,|\,\det \ts{F}>0\}$ and $M^{3\times
3}$ denotes the space of real 3 by 3 matrices. Since we are
interested in modelling large deformations, we assume that the
stored energy function $W$ satisfies that $W\rightarrow\infty$ as
either $\det \ts{F}\rightarrow 0^+$ or
$\norm{\ts{F}}\rightarrow\infty$.

We consider the problem of determining the equilibrium
configuration of the body that satisfies~(\ref{eqn2.3}) a.e.,
and satisfying the boundary condition:
\begin{equation}\label{bcond}
\ts{u}(\ts{x})=\lambda\ts{x},\quad\ts{x}\in\partial\B,
\end{equation}
where $\lambda>0$ is given.

We assume that the stored energy function, in units of energy per unit volume, 
describing the mechanical response of 
the body is given by
\begin{equation}\label{eqn2.5}
W(\ts{x},\ts{F})=\Phi(\ts{x},v_1,v_2,v_3),\quad\ts{F}\in M_+^{3\times 3},\quad 
\ts{x}\in\B,
\end{equation}
for some function $\Phi:\B\times\Real_+^3\To\Real_+$ symmetric in its
last three arguments, and where $v_1,v_2,v_3$ are the eigenvalues of
$(\ts{F}^t\ts{F})^{1/2}$ known as the \textit{principal stretches}. Note that 
for any fixed $\ts{x}$, the material response $W(\ts{x},\cdot)$ corresponds 
to an isotropic and frame indifferent material.

We now restrict attention to the special case in which the deformation
$\ts{u}(\cdot)$ is \textit{radially symmetric}, so that
\begin{equation}\label{eqn2.8}
\ts{u}(\ts{x})=r(R)\,\frac{\ts{x}}{R},\quad\ts{x}\in \B,
\end{equation}
for some scalar function $r(\cdot)$, where $R=\norm{\ts{x}}$. In this
case one can easily check that
\begin{equation}\label{eqn2.9}
v_1=\deru{r}(R),\quad v_2=v_3=\frac{r(R)}{R}.
\end{equation}
Assuming that the dependence of $\Phi$ on $\ts{x}$ in \eqref{eqn2.5} is only on 
$R=\norm{\ts{x}}$, the total stored energy functional, due to internal 
mechanical and gravitational forces (see \cite{CaLe2012}), is given by (up to a 
multiplicative constant of $4\pi$):
\begin{equation}\label{eqn2.10}
I(r)=I_{\mbox{mec}}(r)-I_{\mbox{pot}}(r),
\end{equation}
where 
\begin{eqnarray}
 I_{\mbox{mec}}(r)&=&
\int_0^1 \Phi\left(R,\deru{r}(R),\frac{r(R)}{R}, 
\frac{r(R)}{R}\right)\,R^2\,\dif R,\label{mecfunc}\\
I_{\mbox{pot}}(r)&=&
\int_0^1 \rho_0(R)\dfrac{M_R}{r(R)}\,R^2\,\dif R,\label{potfunc}
\end{eqnarray}
are the mechanical and potential energy functionals respectively. Here 
$\rho_0$ is the mass density of the body (mass per unit volume) in the 
reference configuration, and
\[
 M_R=4\pi\int_0^R\rho_0(u)u^2\,\dif u,
\]
is the mass of the ball in the reference configuration of radius $R$ and 
centered at the origin. We assume that
\begin{equation}\label{massrefden_bdds}
 k_0\le\rho_0(R)\le k_1,\quad 0\le R\le1,
\end{equation}
for some positive constants $k_0$ and $k_1$.

In accord with (\ref{eqn2.3}) we have the inequalities
\begin{equation}\label{eqn2.11}
\deru{r}(R),\,\frac{r(R)}{R}>0,\quad 0<R<1,
\end{equation}
and \eqref{bcond} reduces to:
\begin{equation}\label{eqn2.12}
r(1)=\lambda.
\end{equation}
Our problem now is to minimize the functional $I(\cdot)$ over the set
\begin{eqnarray}
\mathcal{A}_\lambda&=&\big\{r\in
W^{1,1}(0,1)\,|\,r(1)=\lambda,\,\deru{r}(R)>0\mbox{~a.e.
for~}R\in(0,1),\nonumber\\
&&\hspace{1.5in}r(0)\ge0,\, I_{\mbox{mec}}(r)<\infty\big\}.\label{alambda}
\end{eqnarray}
Note that $\mc{A}_\lambda\ne\emptyset$ as $r_\lambda\in\mc{A}_\lambda$ where 
$r_\lambda(R)=\lambda R$.
\section{Existence of minimizers}\label{sec:3}
In this section we show that the functional $I(\cdot)$ in \eqref{eqn2.10} has a 
minimizer over the set $\mc{A}_\lambda$ in \eqref{alambda}. The proofs of the 
results in this section are adaptations of the corresponding ones in 
\cite{Ba82} due to the presence of the potential energy 
functional \eqref{potfunc}. We do emphasize that they are not direct 
consequence of those in \cite{CaLe2012} as these are for maps in Sobolev spaces 
$W^{1,p}$ with $p>3$ and thus they represent continuous deformations.

Throughout this section and the rest of the paper we assume that the stored 
energy function $\Phi$ in \eqref{mecfunc} satisfies that
\begin{equation}\label{hdet_growth}
 \Phi(R,v_1,v_2,v_3)\ge \phi(v_1)+\phi(v_2)+\phi(v_3)+h(v_1v_2v_3),\quad 
R\in[0,1],
\end{equation}
where $\phi,h\,:\,(0,\infty)\To(0,\infty)$ are strictly convex and such that
\begin{enumerate}
\item[H1:]
$\phi(v)\ge C v^\gamma$ for some positive constant $C$ and $1<\gamma<3$;
 \item[H2:]
 $
  \dfrac{h(d)}{d}\To\infty,\quad\mbox{as  }d\To\infty;
 $
\item[H3:]
$
 h(d)\ge Kd^{-s},\quad d>0,
$
for some positive constant $K$ and $s\ge\gamma^*=\frac{\gamma}{\gamma-1}$.
\end{enumerate}
If we let
\[
 \delta_r(R)=\deru{r}(R)\left(\frac{r(R)}{R}\right)^2,
\]
then the specialization of \cite[Eqn. (31)]{CaLe2012} to the radial map 
\eqref{eqn2.8} together with \eqref{massrefden_bdds} gives that
\begin{equation}\label{potfunbdd}
 \abs{\int_0^1 \rho_0(R)\dfrac{M_R}{r(R)}\,R^2\,\dif R}\le
 C\left(\int_0^1\delta_r(R)^{-s}R^2\,\dif R\right)^{\frac{1}{3s}},
\end{equation}
for some positive constant $C$ independent of $r\in\mc{A}_\lambda$. Using this 
we now have the following:
\begin{lem}\label{Ibelow}
Under the growth assumption \eqref{hdet_growth} with H3, the functional 
$I(\cdot)$ 
is bounded below on  $\mc{A}_\lambda$.
\end{lem}
\begin{prueba}
Combining \eqref{massrefden_bdds}, \eqref{hdet_growth} with H3, and 
\eqref{potfunbdd} we 
get for some positive constants $K_1,K_2$ that
\[
 I(r)\ge K_1\int_0^1\delta_r(R)^{-s}R^2\,\dif R-K_2 
\left(\int_0^1\delta_r(R)^{-s}R^2\,\dif R\right)^{\frac{1}{3s}},
\]
for all $r\in\mc{A}_\lambda$. Since the function 
$g(x)=K_1x-K_2x^{\frac{1}{3s}}$ is bounded below for $x\ge0$, the result 
follows.
\end{prueba}
Using this we can now establish the existence of minimizers for $I$ over 
$\mc{A}_\lambda$.
\begin{thm}
 Let the stored energy function $\Phi$ in \eqref{mecfunc} satisfy 
\eqref{hdet_growth} with H1--H3. Then there exists  
$r_\lambda\in\mc{A}_\lambda$ such that
\[
 I(r_\lambda)=\inf_{r\in\mc{A}_\lambda}I(r).
\]
\end{thm}
\begin{prueba}
 Since $\mc{A}_\lambda\ne\emptyset$, it follows from Lemma \ref{Ibelow} that
 $\inf_{r\in\mc{A}_\lambda}I(r)\in\Real$. Let $(r_j)$ with 
$r_j\in\mc{A}_\lambda$ for all $j$, be an infimizing sequence, i.e.,
\[
\inf_{r\in\mc{A}_\lambda}I(r)=\limite{j}{\infty} I(r_j).
\]
Since $(I(r_j))$ is bounded, it follows from the proof of Lemma \ref{Ibelow} 
that the sequence
\begin{equation}\label{lbound_seq}
 \left(K_1\int_0^1\delta_{r_j}(R)^{-s}R^2\,\dif R-K_2 
\left(\int_0^1\delta_{r_j}(R)^{-s}R^2\,\dif R\right)^{\frac{1}{3s}}\right),
\end{equation}
is bounded. Hence the sequence
\[
 \left(\int_0^1\delta_{r_j}(R)^{-s}R^2\,\dif R\right),
\]
must be bounded as well, and thus from \eqref{potfunbdd} that 
$(I_{\mbox{pot}}(r_j))$ 
is bounded. From this and the boundedness of $(I(r_j))$, we get that 
$(I_{\mbox{mec}}(r_j))$ is bounded.

From the boundedness of $(I_{\mbox{mec}}(r_j))$ and \eqref{hdet_growth}, we get 
that
\[
 \left(\int_0^1h(\delta_{r_j}(R))R^2\,\dif R\right),
\]
is bounded. Let $\rho=R^3$ and $u_j(\rho)=r_j^3(\rho^{1/3})$. It follows now 
that
\[
 \dot{u}_j(\rho)=\td{u_j}{\rho}(\rho)=\delta_{r_j}(\rho^{1/3}),
\]
and that the sequence
\[
 \left(\int_0^1h(\dot{u}_j(\rho))\,\dif \rho\right),
\]
is bounded. It follows now from H1 and De La Vall\'ee--Poussin Criterion that 
for some subsequence $(\dot{u}_k)$ of $(\dot{u}_j)$, we have 
$\dot{u}_k\rightharpoonup w$ in $L^1(0,1)$ for some $w\in L^1(0,1)$, and that 
$(\dot{u}_j)$ is equi--integrable. Using H3 is easy to show that $w>0$ 
a.e. Letting
\[
 u(\rho)=\lambda^3-\int_\rho^1w(s)\,\dif s,
\]
we get from the equi--integrability of $(\dot{u}_j)$ that $u_k\To u$ in 
$C[0,1]$. Thus $r_k\To r_\lambda$ in $C[0,1]$ where 
$r_\lambda(R)=u(R^3)^{1/3}$. From these we can conclude $r_k\rightharpoonup 
r_\lambda$ in $W^{1,1}(\eps,1)$ and that 
$\delta_{r_j}\rightharpoonup\delta_{r_\lambda}$ in $L^1(\eps,1)$ for any 
$\eps\in(0,1)$. By the weak lower semi--continuity properties of 
$I_{\mbox{mec}}(\cdot)$ (cf. \cite{BaCuOl1981}), we get that
\begin{eqnarray*}
 \int_\eps^1\!\! \Phi\left[R,\deru{r_\lambda}(R),\frac{r_\lambda(R)}{R}, 
\frac{r_\lambda(R)}{R}\right]R^2\,\dif R&\le&\liminf_k \int_\eps^1\!\! 
\Phi\left[R,\deru{r_k}(R),\frac{r_k(R)}{R}, 
\frac{r_k(R)}{R}\right]R^2\,\dif R,\\ &\le& \liminf_k 
I_{\mbox{mec}}(r_k)<\infty.
\end{eqnarray*}
We get now from the Monotone Convergence Theorem and the arbitrariness of 
$\eps$ that
\begin{equation}\label{wlscmec}
 I_{\mbox{mec}}(r_\lambda)\le \liminf_k 
I_{\mbox{mec}}(r_k).
\end{equation}
This together with the facts that $r_\lambda(0)\ge0$, $\deru{r_\lambda}(R)\ge0$ 
a.e., and $r_\lambda(1)=\lambda$, show that $r_\lambda\in\mc{A}_\lambda$.

To get that $r_\lambda$ is a minimizer of $I$ over $\mc{A}_\lambda$, we must 
still have to deal with the potentials $(I_{\mbox{pot}}(r_k))$. First note that
\begin{eqnarray}
 \abs{I_{\mbox{pot}}(r_k)-I_{\mbox{pot}(r_\lambda)}}=\abs{ 
\int_0^1\dfrac{\rho_0(R)M_RR^2}{r_k(R)r_\lambda(R)}(r_\lambda(R)-r_k(R))\,\dif 
R},\nonumber\\
\le\norm{r_\lambda-r_k}_{C[0,1]}
\left[\int_0^1\dfrac{\rho_0(R)M_RR^2}{r_k^2(R)}\dif R\right]^{\frac{1}{2}}
\left[\int_0^1\dfrac{\rho_0(R)M_RR^2}{r_\lambda^2(R)}\dif 
R\right]^{\frac{1}{2}},\label{potconvg}
\end{eqnarray}
where in the last step we used a weighted Holder's inequality with weight 
$\rho_0(R)M_RR^2$. We now show that each of the two integrals on the right hand 
side of this 
inequality are bounded. Upon recalling \eqref{massrefden_bdds}, we can take 
$M_R\le CR^2$ for some constant $C$. Hence 
\begin{eqnarray*}
 \int_0^1\dfrac{\rho_0(R)M_RR^2}{r_k^2(R)}\dif R&\le& \mbox{(const)} 
\int_0^1\dfrac{R^4}{r_k^2(R)}\,\dif R= \mbox{(const)}
\int_0^1\dfrac{R^2}{\delta_{r_k}(R)}\,\deru{r_k}(R)\dif R,\\
&\le& \mbox{(const)}\left[\int_0^1\dfrac{R^2}{\delta_{r_k}(R)^{\gamma^*}}\dif R 
\right]^{\frac{1}{\gamma^*}}\left[\int_0^1R^2(\deru{r_k}(R))^\gamma\dif R 
\right]^{\frac{1}{\gamma}}.
\end{eqnarray*}
However by the weighted Holder's inequality (with weight $R^2$),
\[
 \left[\int_0^1\dfrac{R^2}{\delta_{r_k}(R)^{\gamma^*}}\dif 
R\right]^{\frac{1}{\gamma^*}}
 \le \mbox{(const)} \left[\int_0^1\dfrac{R^2}{\delta_{r_k}(R)^{s}}\dif 
R\right]^{\frac{1}{s}},
\]
with the sequence on the right hand side bounded. From \eqref{hdet_growth} and 
H1 it follows that
\[
 \int_0^1R^2(\deru{r_k}(R))^\gamma\dif R\le 
\int_0^1R^2\phi(\deru{r_k}(R))\,\dif R\le I_{\mbox{mec}}(r_k),
\]
with the sequence on the right hand side bounded. Combining these results we 
can conclude that the sequence
\[
 \left(\int_0^1\dfrac{\rho_0(R)M_RR^2}{r_k^2(R)}\dif R\right),
\]
is bounded. Combining this with \eqref{potconvg} we get that
\[
 \abs{I_{\mbox{pot}}(r_k)-I_{\mbox{pot}(r_\lambda)}}\le\mbox{(const)} 
\norm{r_\lambda-r_k}_{C[0,1]}\To0,
\]
as $k\To\infty$, which together with \eqref{wlscmec} imply that
\[
 I(r_\lambda)\le \liminf_k I(r_k)=\inf_{r\in\mc{A}_\lambda}I(r),
\]
i.e., that $r_\lambda$ is a minimizer.
\end{prueba}

For our next result we shall need the following assumption: there exist 
constants $M, \eps_0\in(0,\infty)$ such that (cf. \cite{Ba2002})
\begin{equation}\label{derPhibdd}
 \abs{\pd{\Phi}{v_k}(R,\alpha_1v_1,\alpha_2v_2,\alpha_3v_3)v_k}\le 
M\left[\Phi(R,v_1,v_2,v_3)+1\right],
\end{equation}
for all $R\in[0,1]$, $k=1,2,3$, and $\abs{\alpha_i-1}<\eps_0$ for $i=1,2,3$.
The techniques in \cite{Ba82} can now be adapted to show the following 
result.
\begin{thm}\label{ELeqn}
 Let $r$ be any minimizer of $I$ over $\mc{A}_\lambda$. Assume that the
function $\Phi$ satisfies \eqref{derPhibdd}. Then $r\in
C^1(0,1]$, $\deru{r}(R)>0$ for all $R\in(0,1]$, $R^{n-1}\Phi_1(R,r(R))$ is
$C^1(0,1]$, and
\begin{equation}\label{eqn2.13}
\td{}{R}\left[R^{2}\Phi_{,1}(R,r(R))\right]=
2R\Phi_{,2}(R,r(R))+ R^2\,\dfrac{\rho_0(R)M_R}{r^2(R)},
\quad0<R<1,
\end{equation}
subject to~(\ref{eqn2.12}) and $r(0)\ge0$, where:
\begin{equation}\label{Phi-notation}
\Phi_{,i}(R,r(R))=\pd{\Phi}{v_i}\left(R,\deru{r}(R),\displaystyle\frac{r(R)}{R},
\frac{r(R)}{R}\right),\quad i=1,2.
\end{equation}
Moreover, if $r(0)>0$, then 
\begin{equation}\label{cavBC}
 \limite{R}{0^+}R^2\Phi_{,1}(R,r(R))=0.
\end{equation}
\end{thm}

The radial component of the Cauchy stress is given by
\begin{equation}\label{Cstress}
 T(R,r(R))=\dfrac{R^2}{r^2(R)}\Phi_{,1}(R,r(R)).
\end{equation}
Using \eqref{eqn2.13} we now get that
\begin{eqnarray}
\td{T}{R}(R,r(R))&=&2\dfrac{R^2}{r^3(R)}\left[\dfrac{r(R)}{R} 
\Phi_{,2}(R,r(R))-\deru{r}(R)\Phi_{,1}(R,r(R))\right]\nonumber\\
&&~~~~~~~~~~~~~~~~+R^2\,\dfrac{\rho_0(R)M_R}{r^4(R)}.\label{CsDE}
\end{eqnarray}
By the Baker--Ericksen inequality, the right hand side of this equation is 
positive whenever $\deru{r}(R)<\frac{r(R)}{R}$. 

The material of the body $\B$ is \textit{homogeneous} if 
$\rho_0(R)$ is constant, 
still denoted by $\rho_0$, and
\begin{equation}\label{homoSE}
 \Phi(R,v_1,v_2,v_3)=\tilde{\Phi}(v_1,v_2,v_3).
\end{equation}
In this case \eqref{eqn2.13} reduces to
\begin{equation}\label{eqn2.13h}
\td{}{R}\left[R^{2}\tilde{\Phi}_{,1}(r(R))\right]=
2R\tilde{\Phi}_{,2}(r(R))+ \frac{4\pi}{3}\, \rho_0^2\,\dfrac{R^5}{r^2(R)},
\quad0<R<1,
\end{equation}
where now
\begin{equation}\label{Phih-notation}
\tilde{\Phi}_{,i}(r(R))=\pd{\tilde{\Phi}}{v_i}\left(\deru{r}(R),
\displaystyle\frac{r(R)}{R},\frac{r(R)}{R}\right),\quad i=1,2.
\end{equation}
The radial component of the Cauchy stress is now given by
\begin{equation}\label{Cstressh}
 \tilde{T}(r(R))=\dfrac{R^2}{r^2(R)}\tilde{\Phi}_{,1}(r(R)).
\end{equation}
and \eqref{CsDE} reduces to
\begin{equation}\label{CsDEconst}
\td{\tilde{T}}{R}(r(R))=2\dfrac{R^2}{r^3(R)}\left[\dfrac{r(R)}{R} 
\tilde{\Phi}_{,2}(r(R))-\deru{r}(R)\tilde{\Phi}_{,1}(r(R))\right]
+\frac{4\pi}{3}\, \rho_0^2\,\dfrac{R^5}{r^4(R)}.
\end{equation}

\section{Cavitation}\label{sec:4}
In this section we show that when $\lambda$ is sufficiently large, the 
minimizer $r$ of $I(\cdot)$ over $\mathcal{A}_{\lambda}$ has to have $r(0)>0$. 
The proof given here of this fact is an adaptation (to the problem with self 
gravity) of the technique used in \cite{Siva1991} to establish a similar fact 
for compressible inhomogeneous materials. 

We assume for some positive constants $c_0$ and $c_1$,
\begin{equation}\label{hombdd}
 c_0\tilde{\Phi}(v_1,v_2,v_3)\le\Phi(R,v_1,v_2,v_3)\le 
c_1\tilde{\Phi}(v_1,v_2,v_3),
\end{equation}
for all $R\in[0,1]$ and where $\tilde{\Phi}$ corresponds to an isotropic and 
frame indifferent material that satisfies \eqref{hdet_growth} with H1--H3.
We further assume that:
\begin{enumerate}
\item[H4:]
For any $\eta>1$,
\[
\dfrac{v^2}{(v^3-1)^2}\,\tilde{\Phi}\left(\frac{1}{v^2},v,v\right)\in 
L^1(\eta,\infty).
\]
\end{enumerate}
\begin{thm}
 Let $r$ be a minimizer of the functional \eqref{eqn2.10} over 
\eqref{alambda} and assume that \eqref{hombdd} holds where $\tilde{\Phi}$ 
satisfies \eqref{hdet_growth} with H1--H4. Then for $\lambda$ sufficiently 
large we must have that $r(0)>0$.
\end{thm}
\begin{prueba}
 We consider an incompressible deformation given by
 \[
  r_{inc}(R)=\sqrt[3]{R^3+\lambda^3-1}.
 \]
It follows that $r_{inc}\in\mathcal{A}_\lambda$. If $r_\lambda$ is any 
minimizer of $I(\cdot)$ over $\mathcal{A}_\lambda$, then
\begin{eqnarray*}
 \Delta I&=&I(r_{inc})-I(r_\lambda)\\
 &\le&c_1k_1\int_0^1\tilde{\Phi}(r_{inc}(R))R^2\,\dif R-c_0k_0 
\int_0^1\tilde{\Phi}(r_{\lambda}(R))R^2\,\dif R\\
&&-\int_0^1 \rho_0(R)\dfrac{M_R}{r_{inc}(R)}\,R^2\,\dif R+
\int_0^1 \rho_0(R)\dfrac{M_R}{r_\lambda(R)}\,R^2\,\dif R.
\end{eqnarray*}
By \cite[Pro. 4.10]{Si86a} we have that for $\lambda_1\le\lambda_2$
\[
 r_{\lambda_1}(R)\le r_{\lambda_2}(R),\quad0\le R\le1.
\]
Hence for $\lambda_0$ fixed, we get that for $\lambda\ge\lambda_0$,
\[
 \int_0^1 \rho_0(R)\dfrac{M_R}{r_\lambda(R)}\,R^2\,\dif R\le
 \int_0^1 \rho_0(R)\dfrac{M_R}{r_{\lambda_0}(R)}\,R^2\,\dif R.
\]
It follows now that
\begin{eqnarray*}
  \Delta I&\le&c_1k_1\int_0^1\tilde{\Phi}(r_{inc}(R))R^2\,\dif R-c_0k_0 
\int_0^1\tilde{\Phi}(r_{\lambda}(R))R^2\,\dif R\\
&&+\int_0^1 \rho_0(R)\dfrac{M_R}{r_{\lambda_0}(R)}\,R^2\,\dif R.
\end{eqnarray*}
If $r_\lambda(0)=0$, it follows (cf. \cite{Ba82}) that
\[
 \int_0^1\tilde{\Phi}(r_{\lambda}(R))R^2\,\dif R\ge 
\int_0^1\tilde{\Phi}(r_{hom}(R))R^2\,\dif R,
\]
where $r_{hom}(R)=\lambda R$. Hence
\begin{eqnarray*}
  \Delta I&\le&c_1k_1\int_0^1\tilde{\Phi}(r_{inc}(R))R^2\,\dif R-c_0k_0 
\int_0^1\tilde{\Phi}(r_{hom}(R))R^2\,\dif R\\
&&+\int_0^1 \rho_0(R)\dfrac{M_R}{r_{\lambda_0}(R)}\,R^2\,\dif R.
\end{eqnarray*}
Since the third integral on the right hand side of this inequality is fixed, 
the result now follows as in \cite{Siva1991} using H2 and H4.
\end{prueba}

\section{No cavitation results}\label{sec:5}
In this section we give conditions under which the minimizer of $I(\cdot)$ over 
$\mc{A}_\lambda$ for $\lambda$ sufficiently small. must satisfy that $r(0)=0$.

We first consider the case of a homogeneous material for which \eqref{homoSE} 
holds. The functional $I$ is now given by
\begin{equation}\label{homoI}
 I(r)=\int_0^1 \left[\tilde{\Phi}\left(\deru{r}(R),\frac{r(R)}{R}, 
\frac{r(R)}{R}\right)-\rho_0(R)\dfrac{M_R}{r(R)}\right]\,R^2\,\dif R
\end{equation}
We denote by $\lambda_c^h$ the critical boundary displacement for the 
cavitation problem considered in \cite{Ba82} with stored energy function 
$\tilde{\Phi}$.
\begin{thm}\label{nocav1}
 Let $r$ be a minimizer of \eqref{homoI} over $\mc{A}_\lambda$. If 
$\lambda<\lambda_c^h$, then $r(0)=0$.
\end{thm}
\begin{prueba}
 To argue by contradiction, assume that $r(0)>0$. For 
$\hat{\lambda}=(\lambda+\lambda_c^h)/2$, we let
\[
 R_0=\inf\set{R\,|\,r(R)=\hat{\lambda}R}.
\]
Since $\lambda<\lambda_c^h$ and $\frac{r(R)}{R}\To\infty$ as $R\To0^+$, we get 
that $R_0$ is well defined and $R_0>0$. We define
\[
 \hat{r}(R)=\left\{\begin{array}{rcl}\hat{\lambda}R&,&R\le R_0,\\ 
r(R)&,&R>R_0.\end{array}\right.
\]
If follows that $\hat{r}\in\mc{A}_\lambda$. Now
\begin{eqnarray*}
 \Delta I&=&I(r)-I(\hat{r})= \int_0^{R_0} 
\left[\tilde{\Phi}\left(\deru{r}(R),\frac{r(R)}{R}, 
\frac{r(R)}{R}\right)- \tilde{\Phi}\left(\hat{\lambda},\hat{\lambda}, 
\hat{\lambda}\right)\right]R^2\,\dif R\\
&&+\int_0^{R_0}\rho_0(R)M_R\left[\dfrac{1}{\hat{r}(R)}
-\dfrac{1}{r(R)}\right]\,R^2\,\dif R,\\
&\ge&\int_0^{R_0} 
\left[\tilde{\Phi}\left(\deru{r}(R),\frac{r(R)}{R}, 
\frac{r(R)}{R}\right)- \tilde{\Phi}\left(\hat{\lambda},\hat{\lambda}, 
\hat{\lambda}\right)\right]R^2\,\dif R,
\end{eqnarray*}
as $\hat{r}(R)\le r(R)$ for $R\le R_0$. Since $\hat{\lambda}<\lambda_c^h$ and 
$r(0)>0$, the results in \cite{Ba82} imply that
\[
 \int_0^{R_0}\tilde{\Phi}\left(\deru{r}(R),\frac{r(R)}{R}, 
\frac{r(R)}{R}\right)R^2\,\dif R>
\int_0^{R_0} 
\tilde{\Phi}\left(\hat{\lambda},\hat{\lambda}, 
\hat{\lambda}\right)R^2\,\dif R,
\]
and thus that $\Delta I>0$ which contradicts the minimality of $r$.
\end{prueba}

For the next result we take $\beta\equiv0$ and $\gamma\equiv1$ in 
\eqref{inhomoSE}. We let $d_0$ be the value of the argument at which $h$ 
assumes its global minimum value. Note that $\deru{h}(d)<0$ for $d<d_0$. 
\begin{thm}
Let $r$ be a minimizer of \eqref{eqn2.10} over $\mc{A}_\lambda$ corresponding 
to the stored energy function \eqref{inhomoSE}. Assume that 
$\deru{\alpha}(R)\le0$ for all $R$ and that $\deru{\phi}(t)\ge0$ for all $t$.
Then if $\lambda^3<d_0$ we must have that $r(0)=0$.
\end{thm}
\begin{prueba}
 As in the proof of Theorem \ref{nocav1}, we argue by contradiction. Thus we 
assume that 
$r(0)>0$ and let $\hat{\lambda}=(\lambda+d_0^{\frac{1}{3}})/2$. Define $R_0$ 
and $\hat{r}$ as in the proof of Theorem \ref{nocav1}. Then
\begin{eqnarray*}
 \Delta I&=&I(r)-I(\hat{r})= \int_0^{R_0} 
\left[\Phi\left(R,\deru{r}(R),\frac{r(R)}{R}, 
\frac{r(R)}{R}\right)- \Phi\left(R,\hat{\lambda},\hat{\lambda}, 
\hat{\lambda}\right)\right]R^2\,\dif R\\
&&+\int_0^{R_0}\rho_0(R)M_R\left[\dfrac{1}{\hat{r}(R)}
-\dfrac{1}{r(R)}\right]\,R^2\,\dif R,\\
&\ge&\int_0^{R_0} 
\left[\Phi\left(R,\deru{r}(R),\frac{r(R)}{R}, 
\frac{r(R)}{R}\right)- \Phi\left(R,\hat{\lambda},\hat{\lambda}, 
\hat{\lambda}\right)\right]R^2\,\dif R,
\end{eqnarray*}
as $\hat{r}(R)\le r(R)$ for $R\le R_0$. Using the convexity of $\phi$, 
$\psi$, and $h$ in \eqref{inhomoSE}, we have now that (see \cite[Page 
589]{Ba82}):
\begin{eqnarray*}
&\left[\Phi\left(R,\deru{r}(R),\frac{r(R)}{R}, 
\frac{r(R)}{R}\right)- \Phi\left(R,\hat{\lambda},\hat{\lambda}, 
\hat{\lambda}\right)\right]R^2\ge\hspace{2.5in}&\\
&\hspace{1.0in}\alpha(R) 
\deru{\phi}(\hat{\lambda})\left(R^2\deru{r}(R)+2Rr(R)-3\hat{\lambda}
R^2\right)+
\deru{h}(\hat{\lambda}^3)\left(r^2(R)\deru{r}(R)-\hat{\lambda}^3R^2\right)
,&\\
&=\alpha(R)\deru{\phi}(\hat{\lambda})\deru{(r(R)R^2-\hat{\lambda}R^3)}
+\frac{1}{3}\deru{h}(\hat{\lambda}^3)\deru{(r^3(R)-\hat{\lambda}^3R^3)}.&
\end{eqnarray*}
It follows now, after integrating by parts, that
\begin{eqnarray*}
 &\displaystyle\int_0^{R_0} 
\left[\Phi\left(R,\deru{r}(R),\frac{r(R)}{R}, 
\frac{r(R)}{R}\right)- \Phi\left(R,\hat{\lambda},\hat{\lambda}, 
\hat{\lambda}\right)\right]R^2\,\dif R\ge \hspace{2.0in}&\\ 
&\displaystyle-\deru{\phi}(\hat{\lambda})\int_0^{R_0}
\deru{\alpha}(R)(r(R)R^2-\hat{\lambda}
R^3)\,\dif R-\frac{1}{3}\deru{h}(\hat{\lambda}^3)r^3(0)\ge
-\frac{1}{3}\deru{h}(\hat{\lambda}^3)r^3(0)>0.&
\end{eqnarray*}
Thus $\Delta I>0$ which contradicts the minimality of $r$.
\end{prueba}

Our next results is for inhomogeneous materials of the form
\begin{equation}\label{inhomoSE}
 \Phi(R,v_1,v_2,v_3)=\alpha(R)\sum_i\phi(v_i)+\beta(R)\sum_{i<j}\psi(v_iv_j)+ 
\gamma(R)h(v_1v_2v_3),
\end{equation}
where $\alpha$, $\beta$, and $\gamma$ are smooth positive functions over 
$[0,1]$, $\phi$ and $\psi$ are nonegative convex functions and with $h$ 
strictly convex. In this case it easy to see that for some constant $C>0$,
\begin{equation}\label{PhiRbdd}
\abs{\pd{\Phi}{R}(R,v_1,v_2,v_3)}\le C\,\tilde{\Phi}(v_1,v_2,v_3),
\end{equation}
where
\begin{equation}\label{homoSE2}
\tilde{\Phi}(v_1,v_2,v_3)=\sum_i\phi(v_i)+\sum_{i<j}\psi(v_iv_j)+h(v_1v_2v_3).
\end{equation}
We now denote by $\lambda_c^h$ the critical boundary displacement 
corresponding to the stored energy function $\Phi(0,v_1,v_2,v_3)$. Note that 
any deformation with finite $\Phi(0,v_1,v_2,v_3)$ energy, has 
finite $\tilde{\Phi}$ energy as well. The following result is reminiscent to 
\cite[Proposition 12]{Siva1991}. It shows that any minimizer which leaves the 
centre intact must have strains at the origin less than $\lambda_c^h$.

\begin{thm}\label{nocav_ls}
 Let $r\in\mc{A}_\lambda$ satisfy $r(0)=0$ and that $\ell\in(0,\infty]$ where
 \[
  \limite{R}{0^+}\deru{r}(R)=\ell.
 \]
If $\ell>\lambda_c^h$, then $r$ can not be a minimizer of $I(\cdot)$ over 
$\mc{A}_\lambda$.
\end{thm}
\begin{prueba}
The following proof is similar to the one of \cite[Proposition 12]{Siva1991} 
except for the treatment of the gravitational potential and the specific 
dependence on $R$ of the stored energy function.

Fro $\eps\in(0,1)$ we let
\[
 \lambda(\eps)=\dfrac{r(\eps)}{\eps}.
\]
From the given hypotheses, it follows that
\[
 \lambda(\eps)\To\ell,\quad\mbox{as }\eps\To0^+.
\]
Assume for the moment that $\ell$ is finite, and let $r_c$ be a cavitating 
extrema corresponding to $\Phi(0,v_1,v_2,v_3)$. We define
\begin{equation}\label{testfunc}
 r_\eps(R)=\left\{\begin{array}{rcl} \alpha_\eps 
r_c\left(\frac{R}{\alpha_\eps}\right)&,&R\in[0,\eps],\\ r(R)&,&R\in(\eps,1],
\end{array}\right.
\end{equation}
where $\alpha_\eps$ is such that $\alpha_\eps 
r_c(\eps/\alpha_\eps)=\lambda(\eps)\eps$. That $\alpha_\eps$ exists follows 
from the fact that $\lambda(\eps)>\lambda_c^h$ for $\eps$ sufficiently small.
Now
\begin{eqnarray*}
 \Delta I&=&I(r_\eps)-I(r)\\
 &=&\int_0^\eps 
R^2\bigg[\Phi\left(R,\deru{r_c}(R/\alpha_\eps),\frac{\alpha_\eps 
r_c(R/\alpha_\eps)}{R},\frac{\alpha_\eps r_c(R/\alpha_\eps)}{R}\right)\\
&&\quad\quad-\Phi\left(R,\deru{r}(R),\frac{r(R)}{R}, 
\frac{r(R)}{R}\right)\bigg]\,\dif R\\
&&-\int_0^\eps\rho_0(R)\dfrac{M_R}{\alpha_\eps r_c(R/\alpha_\eps)} 
\,R^2\,\dif R
+\int_0^\eps\rho_0(R)\dfrac{M_R}{r(R)}\,R^2\,\dif R
\end{eqnarray*}
With the change of variables $R=\eps U$ with $U\in[0,1]$, we can write the 
above as
\begin{eqnarray*}
 \Delta I&=&\eps^3\int_0^1 
U^2\bigg[\Phi\left(\eps U,\deru{r_c}(\eps U/\alpha_\eps),\frac{\alpha_\eps 
r_c(\eps U/\alpha_\eps)}{\eps U},\frac{\alpha_\eps 
r_c(\eps U/\alpha_\eps)}{\eps U}\right)\\
&&\quad\quad-\Phi\left(\eps U,\deru{r}(\eps U),\frac{r(\eps U)}{\eps U}, 
\frac{r(\eps U)}{\eps U}\right)\bigg]\,\dif U\\
&&-\int_0^\eps\rho_0(R)\dfrac{M_R}{\alpha_\eps r_c(R/\alpha_\eps)} 
\,R^2\,\dif R
+\int_0^\eps\rho_0(R)\dfrac{M_R}{r(R)}\,R^2\,\dif R
\end{eqnarray*}
We first examine the gravitational integrals. For this we use that 
$\rho_0(\cdot)$ is nonnegative and bounded above, and that $M_R$ is bounded by 
a constant times $R^3$. Since
\[
 \dfrac{r_c(\eps/\alpha_\eps)}{\eps/\alpha_\eps}=\lambda(\eps)
 \To\ell\quad\mbox{as }\eps\To0^+,
\]
we get that
\[
 \dfrac{\eps}{\alpha_\eps}\To\mu,
\]
where $\mu>0$ and $r_c(\mu)/\mu=\ell$. Upon recalling that $\frac{r_c(S)}{S}$ 
is a decreasing function of $S$, we have that:
\begin{eqnarray*}
 \int_0^\eps\rho_0(R)\dfrac{M_R}{\alpha_\eps r_c(R/\alpha_\eps)} 
\,R^2\,\dif R&=& \eps^3\int_0^1\rho_0(\eps U)\dfrac{M_{\eps U}}{\alpha_\eps 
r_c(\eps U/\alpha_\eps)} 
\,U^2\,\dif U\\
&\le& K_1\eps^5\int_0^1\dfrac{U^4}{\alpha_\eps r_c(\eps/\alpha_\eps)/\eps} 
\,\dif U\le \dfrac{K\eps^5}{5\ell}.
\end{eqnarray*}
As $r(0)=0$ and $\ell>0$, the function $r(R)/R$ is positive and continuous in 
$[0,1]$. Thus if $v_0$ is its minimum value, we have that for some positive 
constant $L$:
\[
 \int_0^\eps\rho_0(R)\dfrac{M_R}{r(R)}\,R^2\,\dif R= \eps^3\int_0^1
 \rho_0(\eps U)\dfrac{M_{\eps U}}{r(\eps U)}\,U^2\,\dif U\le 
\frac{M\eps^5}{5v_0}.
\]
Thus both gravitational potential terms go to zero faster than $\eps^3$.

We now examine the mechanical potential terms in $\Delta I$. For this we note 
that
\begin{eqnarray}
 &\displaystyle\int_0^1 
U^2\bigg[\Phi\left(\eps U,\deru{r_c}(\eps U/\alpha_\eps),\frac{\alpha_\eps 
r_c(\eps U/\alpha_\eps)}{\eps U},\frac{\alpha_\eps 
r_c(\eps U/\alpha_\eps)}{\eps U}\right)\hspace{2.0in}&\nonumber\\
&\quad\quad\quad\quad\hspace{1.0in}
-\Phi\left(\eps U,\deru{r}(\eps U),\frac{r(\eps U)}{\eps 
U}, 
\frac{r(\eps U)}{\eps U}\right)\bigg]\,\dif U=&\nonumber\\
&\displaystyle\int_0^1 
U^2\bigg\{\bigg[\Phi\left(\eps U,\deru{r_c}(\eps 
U/\alpha_\eps),\frac{\alpha_\eps 
r_c(\eps U/\alpha_\eps)}{\eps U},\frac{\alpha_\eps 
r_c(\eps U/\alpha_\eps)}{\eps U}\right)\hspace{2.0in}&\nonumber\\
&\displaystyle-\Phi\left(0,\deru{r_c}(\eps U/\alpha_\eps),\frac{\alpha_\eps 
r_c(\eps U/\alpha_\eps)}{\eps U},\frac{\alpha_\eps 
r_c(\eps U/\alpha_\eps)}{\eps U}\right)\bigg]&\nonumber\\
&\displaystyle+\bigg[\Phi\left(0,\deru{r_c}(\eps 
U/\alpha_\eps),\frac{\alpha_\eps 
r_c(\eps U/\alpha_\eps)}{\eps U},\frac{\alpha_\eps 
r_c(\eps U/\alpha_\eps)}{\eps U}\right)\hspace{2.0in}&\nonumber\\
&\quad\quad\quad\quad-\Phi(0,\lambda(\eps),\lambda(\eps),\lambda(\eps))\bigg]
&\nonumber\\
&+\bigg[\Phi(0,\lambda(\eps),\lambda(\eps),\lambda(\eps))
\displaystyle-\Phi\left(\eps U,\deru{r}(\eps U),\frac{r(\eps U)}{\eps U}, 
\frac{r(\eps U)}{\eps U}\right)\bigg]\bigg\}\,\dif 
U\hspace{1.0in}&\label{mecpotM}
\end{eqnarray}
From \eqref{PhiRbdd} and Taylor's Theorem, we have that
\begin{eqnarray*}
 &\displaystyle\int_0^1 
U^2\bigg[\Phi\left(\eps U,\deru{r_c}(\eps 
U/\alpha_\eps),\frac{\alpha_\eps 
r_c(\eps U/\alpha_\eps)}{\eps U},\frac{\alpha_\eps 
r_c(\eps U/\alpha_\eps)}{\eps U}\right)\hspace{2.0in}&\\
&\displaystyle-\Phi\left(0,\deru{r_c}(\eps U/\alpha_\eps),\frac{\alpha_\eps 
r_c(\eps U/\alpha_\eps)}{\eps U},\frac{\alpha_\eps 
r_c(\eps U/\alpha_\eps)}{\eps U}\right)\bigg]\,\dif U&\\
&\le C_1\eps \displaystyle\int_0^1 
U^3\tilde{\Phi}\left(\deru{r_c}(\eps 
U/\alpha_\eps),\frac{\alpha_\eps 
r_c(\eps U/\alpha_\eps)}{\eps U},\frac{\alpha_\eps 
r_c(\eps U/\alpha_\eps)}{\eps U}\right)\,\dif U\le C_2\eps,&
\end{eqnarray*}
where we used that $r_c$ has finite $\tilde{\Phi}$ energy. Thus the first 
bracketed term in \eqref{mecpotM} goes to zero with $\eps$. For the third term, 
we note that the functions $\deru{r}(S)$ and $\frac{r(S)}{S}$ are $C[0,1]$ and 
positive. Hence for some $M>0$,
\[
 \abs{\Phi\left(\eps U,\deru{r}(\eps U),\frac{r(\eps U)}{\eps U}, 
\frac{r(\eps U)}{\eps U}\right)}\le M,\quad\forall U,
\]
and since
\[
 \Phi\left(\eps U,\deru{r}(\eps U),\frac{r(\eps U)}{\eps U}, 
\frac{r(\eps U)}{\eps U}\right)\To\Phi(0,\ell,\ell,\ell),
\]
pointwise, we get by the Lebesgue dominated convergence theorem that
\[
 \int_0^1U^2\Phi\left(\eps U,\deru{r}(\eps U),\frac{r(\eps U)}{\eps U}, 
\frac{r(\eps U)}{\eps U}\right)\,\dif U\To 
\int_0^1U^2\Phi(0,\ell,\ell,\ell)\,\dif U,
\]
as $\eps\To0^+$. This together with $\lambda(\eps)\To\ell$ yields that
\[
\int_0^1U^2\bigg[\Phi(0,\lambda(\eps),\lambda(\eps),\lambda(\eps))
\displaystyle-\Phi\left(\eps U,\deru{r}(\eps U),\frac{r(\eps U)}{\eps U}, 
\frac{r(\eps U)}{\eps U}\right)\bigg]\,\dif U\To0,
\]
as $\eps\To0^+$. Thus the third term in \eqref{mecpotM} goes to zero with 
$\eps$ as well.

For the second term in \eqref{mecpotM}, note that with the change of variables 
$Z=(\eps/\alpha_\eps)U$, we get
\begin{eqnarray*}
 &\displaystyle\int_0^1U^2\,\Phi\left(0,\deru{r_c}(\eps 
U/\alpha_\eps),\frac{\alpha_\eps 
r_c(\eps U/\alpha_\eps)}{\eps U},\frac{\alpha_\eps 
r_c(\eps U/\alpha_\eps)}{\eps U}\right)\,\dif U=\hspace{1.5in}&\\
&\displaystyle\left(\frac{\alpha_\eps}{\eps}\right)^3
\int_0^{\eps/\alpha_\eps}Z^2\,\Phi\left(0,\deru{r_c}(Z),\frac{
r_c(Z)}{Z},\frac{r_c(Z)}{Z}\right)\,\dif Z&\\
&\displaystyle\To\frac{1}{\mu^3}
\int_0^{\mu}Z^2\,\Phi\left(0,\deru{r_c}(Z),\frac{
r_c(Z)}{Z},\frac{r_c(Z)}{Z}\right)\,\dif Z&\\
&\displaystyle=\int_0^1U^2\,\Phi\left(0,\deru{r_c}(\mu U) ,\frac{ 
r_c(\mu U)}{\mu U},\frac{r_c(\mu U)}{\mu 
U}\right)\,\dif U.
\end{eqnarray*}
It follows now that
\begin{eqnarray*}
&\displaystyle\int_0^1U^2\bigg[\Phi\left(0,\deru{r_c}(\eps 
U/\alpha_\eps),\frac{\alpha_\eps 
r_c(\eps U/\alpha_\eps)}{\eps U},\frac{\alpha_\eps 
r_c(\eps U/\alpha_\eps)}{\eps U}\right)\hspace{2.0in}&\\
&\quad\quad\quad\quad-\Phi(0,\lambda(\eps),\lambda(\eps),\lambda(\eps))\bigg]
\,\dif R\To&\\
&\displaystyle\int_0^1U^2\bigg[\Phi\left(0,\deru{r_c}(\mu U) ,\frac{ 
r_c(\mu U)}{\mu U},\frac{r_c(\mu U)}{\mu 
U}\right)-\Phi(0,\ell,\ell,\ell)\bigg]\,\dif R<0,&
\end{eqnarray*}
where the last inequality follows since 
$\ell>\lambda_c^h$ and $\tilde{r}(U)=\mu^{-1}r_c(\mu U)$ is the minimizer for 
the functional with stored energy $\Phi(0,v_1,v_2,v_3)$ and boundary condition 
$\tilde{r}(1)=\ell$. Collecting all of the intermediate results so far, we get 
that $\eps^{-3}\Delta I<0$ for $\eps$ sufficiently small, which contradicts the 
minimality of $r$.

The case $\ell=\infty$, can be handled in a similar fashion using a 
suitable incompressible deformation on $[0,\eps]$ in \eqref{testfunc}. See 
\cite[Proposition 12]{Siva1991} for details.
\end{prueba}

\section{Numerical results}\label{sec:6}
In this section we present some numerical results that confirm some of the 
results of previous sections. We employ two numerical schemes: a descent method 
for the minimization of \eqref{eqn2.10} based on a gradient flow iteration; and 
a shooting method that solves directly the Euler-Lagrange boundary value 
problem \eqref{eqn2.12}, \eqref{eqn2.13}, and \eqref{cavBC}. The use of 
adaptive ode solvers in the shooting method allows for a more precise 
computation of the equilibrium states, especially near $R=0$ where both strains 
in our problem tend to develop boundary layers. After the equilibrium is 
computed via the shooting method, it is compared to the results of the 
descent iteration in order to get some assurance of its minimizing character.

A \textit{gradient flow} iteration (cf. \cite{Neu1997}) assumes that $r$ 
depends on a flow parameter 
$t$, and that $r(R,t)$ satisfies 
\begin{eqnarray}
\displaystyle\tdd{}{R}(r_t(R,t))&=&-\td{}{R}\left[R^{2}\Phi_{,1}
(R,r(R,t))\right]
+2R\Phi_{,2}(R,r(R,t))\label{GFde}\\
&&~~~~~~~~~~~+R^2\,\dfrac{\rho_0(R)M_R}{r^2(R,t)},\quad0<R<1,~t>0,\nonumber\\
r(1,t)&=&\lambda,\quad 
\limite{R}{0^+}\left[\td{}{R}(r_t(R,t))+R^2\Phi_{,1}(R,r(R,t))\right]=0, 
\quad t\ge0.\label{GFbc}
\end{eqnarray}
(Here $r_t=\pd{r}{t}$.) The gradient flow equation leads to a descent method for
the minimization of \eqref{eqn2.10} over \eqref{alambda}. After 
discretization of the partial derivative with respect to ``$t$", one can use a
finite element method to solve the resulting flow
equation. In particular, if we let $\Delta t>0$
be given, and set $t_{i+1}=t_i+\Delta t$ where $t_0=0$, we can approximate
$r_t(R,t_i)$ with:
\[
z_i(R)= \dfrac{r_{i+1}(R)-r_i(R)}{\Delta
t},
\]
where $r_i(R)=r(R,t_i)$, etc. (We take $r_0(R)$
to be some initial deformation satisfying the boundary condition at $R=1$
, e.g., $\lambda R$.) Inserting this approximation into the weak form
of
\eqref{GFde}, \eqref{GFbc}, and evaluating the right hand side of
\eqref{GFde} at $r=r_i$, we arrive at the following iterative
formula:
\begin{eqnarray}
&\displaystyle\int_0^1\deru{z}_i(R)\deru{v}(R)\,\dif R+ 
\int_0^1\Big[R^{2}\Phi_{,1}(R,r_i(R))\deru{v}(R)
~~~~~~~~~~~~~~~~~~~~~~~~~~~~~~~~~\nonumber\\
&+\left(2R\Phi_{,2}(R,r_i(R))+
 R^2\,\dfrac{\rho_0(R)M_R}{r_i^2(R)}\right)v(R)\Big]\,\dif R=0,
\label{floweq1}&
\end{eqnarray}
for all functions $v$ such that $v(1)=0$ and sufficiently smooth so
that the
integrals above are well defined. Given $r_i$, one can solve the above
equation for $z_i$ via some
finite element scheme, and then set $r_{i+1}=r_i+\Delta t\,z_i$. This 
process is repeated 
for $i=0,1,\ldots$, until
$r_{i+1}-r_i$ is ``small'' enough, or some maximum value of
``$t$'' is reached, declaring the last $r_i$ as an approximate minimizer of 
\eqref{eqn2.10} over \eqref{alambda}.

In the shooting method technique, for given $\nu>0$, we solve the initial value 
problem
\begin{subequations}\label{shootIVP}
\begin{eqnarray}
 &\displaystyle\td{}{R}\left[R^{2}\Phi_{,1}(R,r(R))\right]=
2R\Phi_{,2}(R,r(R))+ R^2\,\dfrac{\rho_0(R)M_R}{r^2(R)},
\quad0<R<1,\label{shootIVP1}&\\
&r(1)=\lambda,\quad\deru{r}(1)=\nu,\label{shootIVP2}&
\end{eqnarray}
\end{subequations}
from $R=1$ to $R=0$. The value of $\nu$ is adjusted so that
\begin{equation}\label{shootIVP3}
\limite{R}{0^+}R^2\Phi_{,1}\left(R,\nu,\dfrac{r(R)}{R},
\dfrac{r(R)}{R}\right)=0.
\end{equation}
In actual calculations we solve \eqref{shootIVP} from $R=1$ to $R=\eps$, where 
$\eps>0$ is small, and replace \eqref{shootIVP3} with
\begin{equation}\label{discT}
 \Phi_{,1}\left(\eps,\nu,\dfrac{r(\eps)}{\eps},
\dfrac{r(\eps)}{\eps}\right)=0.
\end{equation}
This equation is solved for $\nu$ via a secant type iteration which requires 
repeated solutions of the initial value problem \eqref{shootIVP} from $R=1$ to 
$R=\eps$. These intermediate initial value problems are solved with the routine 
\texttt{ode45} of the MATLAB\texttrademark~ ode suite.

Our first set of simulations are for the homogeneous case \eqref{homoI} and for 
which the stored energy function $\tilde{\Phi}$ is given by
\[
\tilde{\Phi}(v_1,v_2, v_3)=\frac{\kappa}{p}\,(v_1^p+v_2^p+v_3^p)+h(v_1v_2v_3),
\]
for which
\[
h(d)=C\,d^{\gamma}+D\,d^{-\delta},
\]
where $p<3$, $C\ge0$, $D\ge0$ and $\gamma,\delta>0$. The reference configuration
is mechanically stress free provided:
\[
D=\dfrac{\kappa+C\gamma}{\delta}.
\]
The mass density function $\rho_0$ is taken to be a constant. For the 
simulations we used the following values of the mechanical parameters in 
$\tilde{\Phi}$:
\[
 p=2,\quad\kappa=1,\quad C=1,\quad\gamma=\delta=2,
\]
with a value of $\eps=0.001$ in \eqref{discT} and in the shooting method. The 
gradient flow iteration was used as a predictor for the shooting method, with 
the integrals in \eqref{floweq1} computed over $(\eps,1)$ as well. 

\begin{figure}
\begin{center}
\scalebox{0.5}{\includegraphics{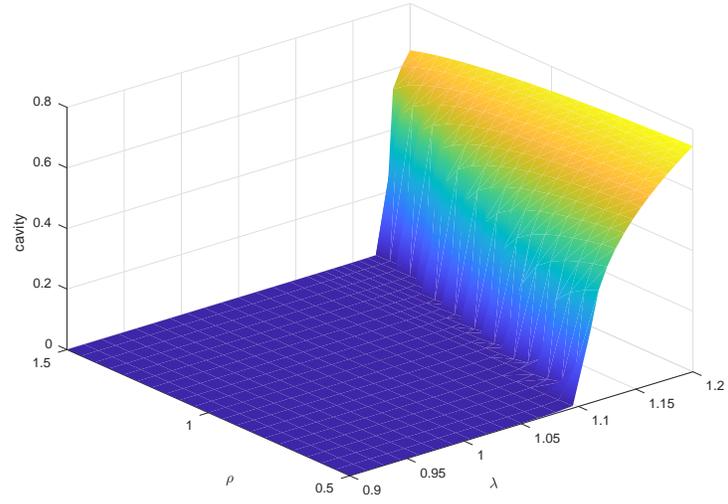}}
\end{center}
\caption{Cavity surface.}
\label{fig:1}
\end{figure}

\begin{figure}
\begin{center}
\scalebox{0.5}{\includegraphics{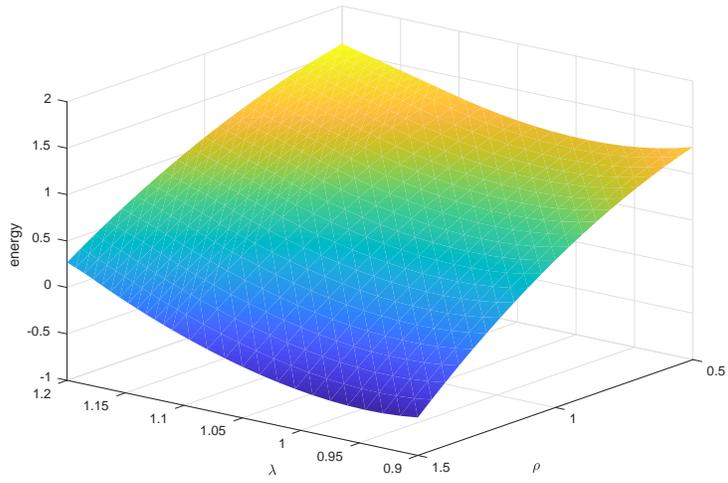}}
\end{center}
\caption{Energy surface.}
\label{fig:2}
\end{figure}

In our first 
simulation we show the approximate cavity radius as a function of $\rho_0$ and 
$\lambda$, for values of $\rho_0\in[0.5,1.5]$ and $\lambda\in[0.9,1.2]$. The 
resulting surface is shown in Figure \ref{fig:1}. We note that for fixed values 
of $\rho_0$, the cavity size is essentially zero up to some certain value of 
$\lambda$ (the critical boundary displacement corresponding to $\rho_0$), after 
which the graph becomes concave. This critical boundary displacement appears to 
be an increasing function of $\rho_0$. In Figure \ref{fig:2} we show the 
corresponding surface for the energies of the approximate minimizers. For fixed 
values of $\lambda$, the energy is an increasing function of $\rho_0$, while 
for $\rho_0$ constant the energy is non monotone with respect to $\lambda$ with 
a convex shape.

\begin{figure}
\begin{center}
\scalebox{0.75}{\includegraphics{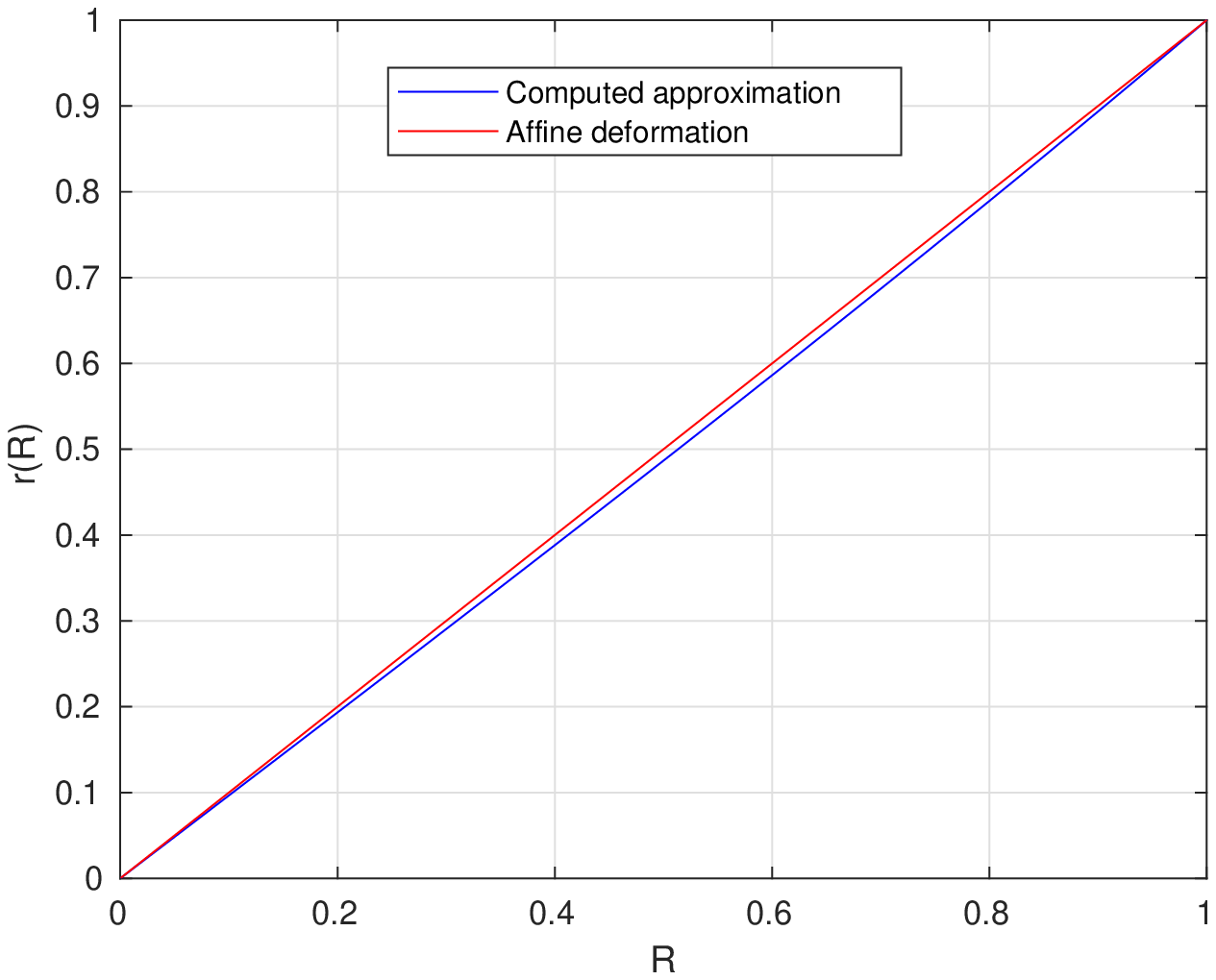}}
\end{center}
\caption{Radial displacement for $\lambda=1$ and $\rho_0=1$.}
\label{fig:3}
\end{figure}

\begin{figure}
\begin{center}
\scalebox{0.75}{\includegraphics{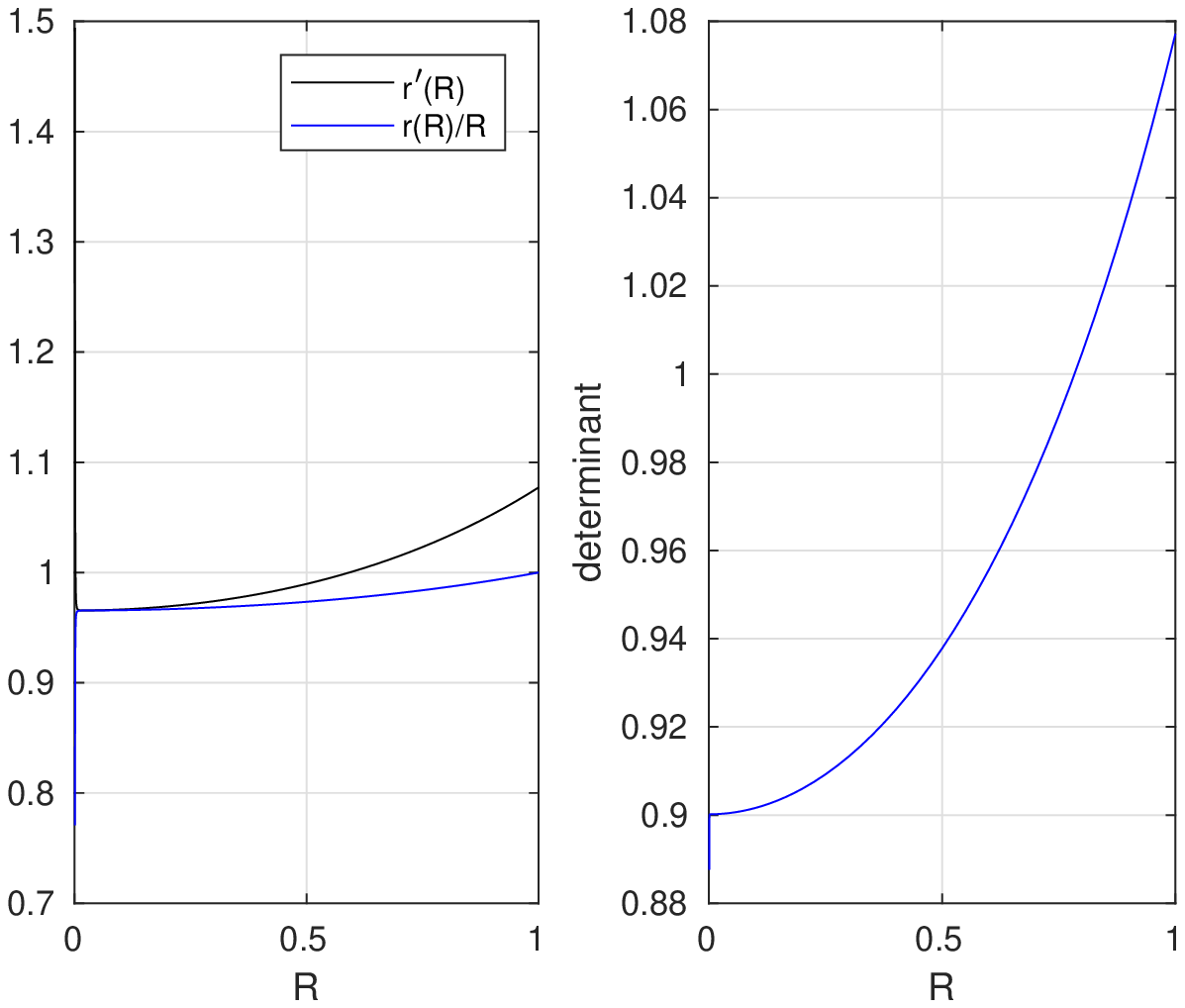}}
\end{center}
\caption{Strains and determinant for $\lambda=1$ and $\rho_0=1$.}
\label{fig:4}
\end{figure}

\begin{figure}
\begin{center}
\scalebox{0.75}{\includegraphics{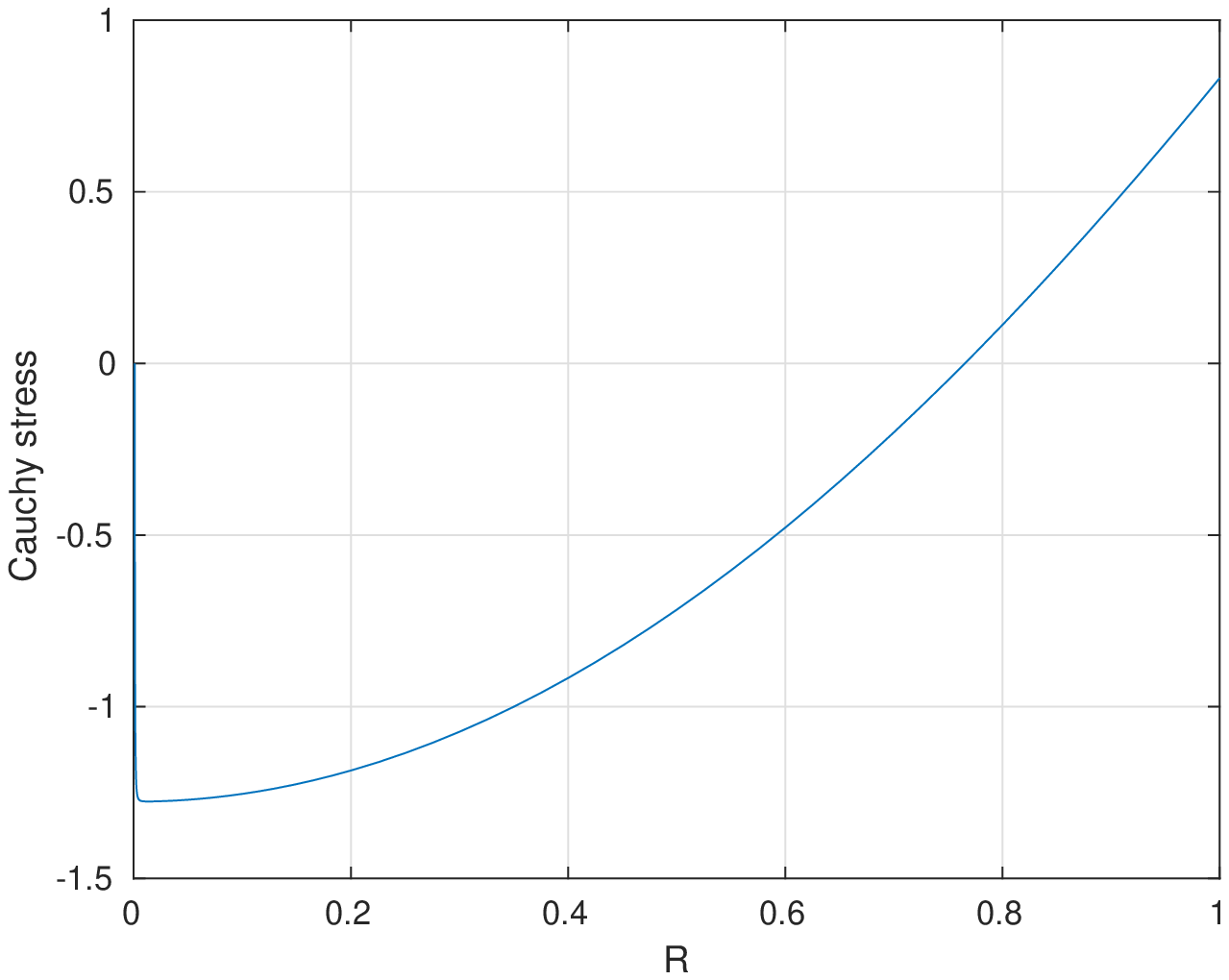}}
\end{center}
\caption{Cauchy stress for $\lambda=1$ and $\rho_0=1$.}
\label{fig:5}
\end{figure}

Our next simulations are for fixed values of $\lambda$ and $\rho_0$. In the 
first case $\lambda=1$ and $\rho_0=1$. In Figure \ref{fig:3} we show the 
computed minimizer $r$ compared to the affine deformation $\lambda R$. The 
value of $r(0.001)$ is $7.7078\times10^{-4}$ with an energy of $0.49074$. Figure 
\ref{fig:4} shows plots of the strains $\deru{r}(R)$ and $\frac{r(R)}{R}$, and 
the determinant $\deru{r}(R)(r(R)/R)^2$ in this case. Also in Figure 
\ref{fig:5} we show the graph of the corresponding Cauchy stress 
\eqref{Cstressh}. We note the boundary layer close to $R=\eps$ in these plots. 
This boundary layer is a numerical artefact since the numerical scheme tries to 
make the Cauchy stress zero, while for this value of $\lambda$ the value of 
$r(0)$ should be zero. Still in this case the numerical scheme converges to the 
solution with $r(0)=0$ as $\eps\To0^+$ (cf. \cite[Section 5]{Si86a}).

\begin{figure}
\begin{center}
\scalebox{0.75}{\includegraphics{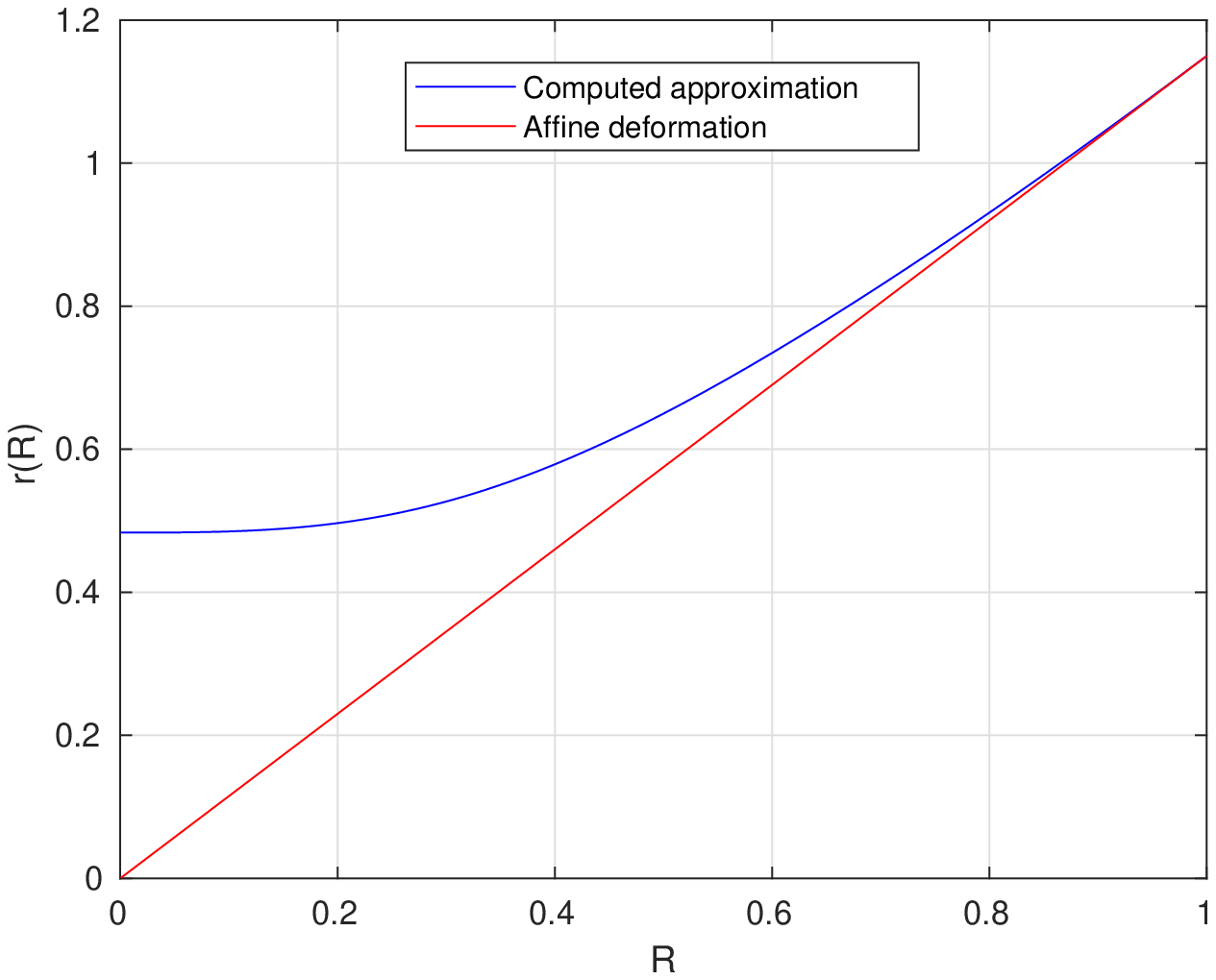}}
\end{center}
\caption{Radial displacement for $\lambda=1.15$ and $\rho_0=1$.}
\label{fig:6}
\end{figure}

\begin{figure}
\begin{center}
\scalebox{0.75}{\includegraphics{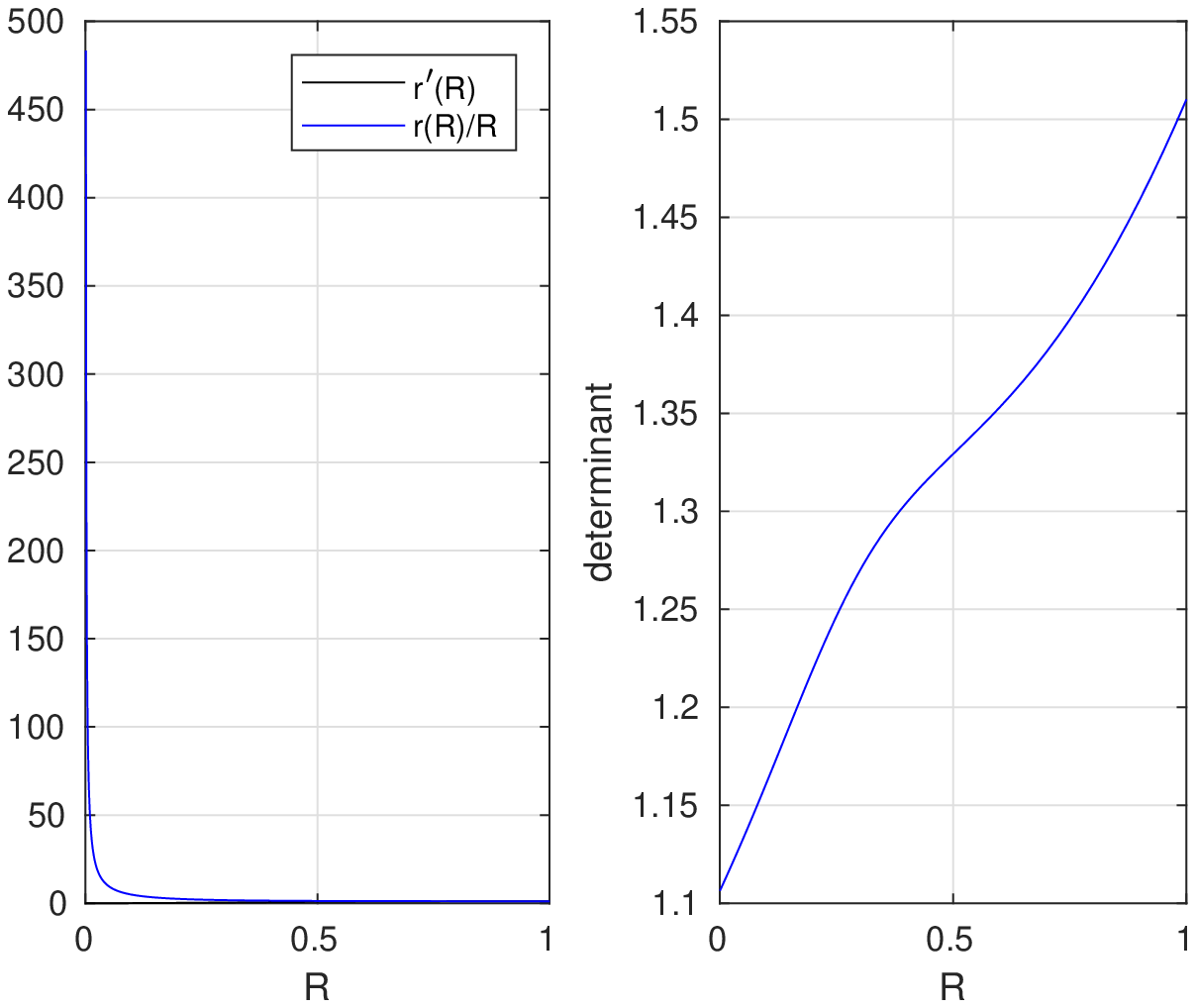}}
\end{center}
\caption{Strains and determinant for $\lambda=1.15$ and $\rho_0=1$.}
\label{fig:7}
\end{figure}

\begin{figure}
\begin{center}
\scalebox{0.75}{\includegraphics{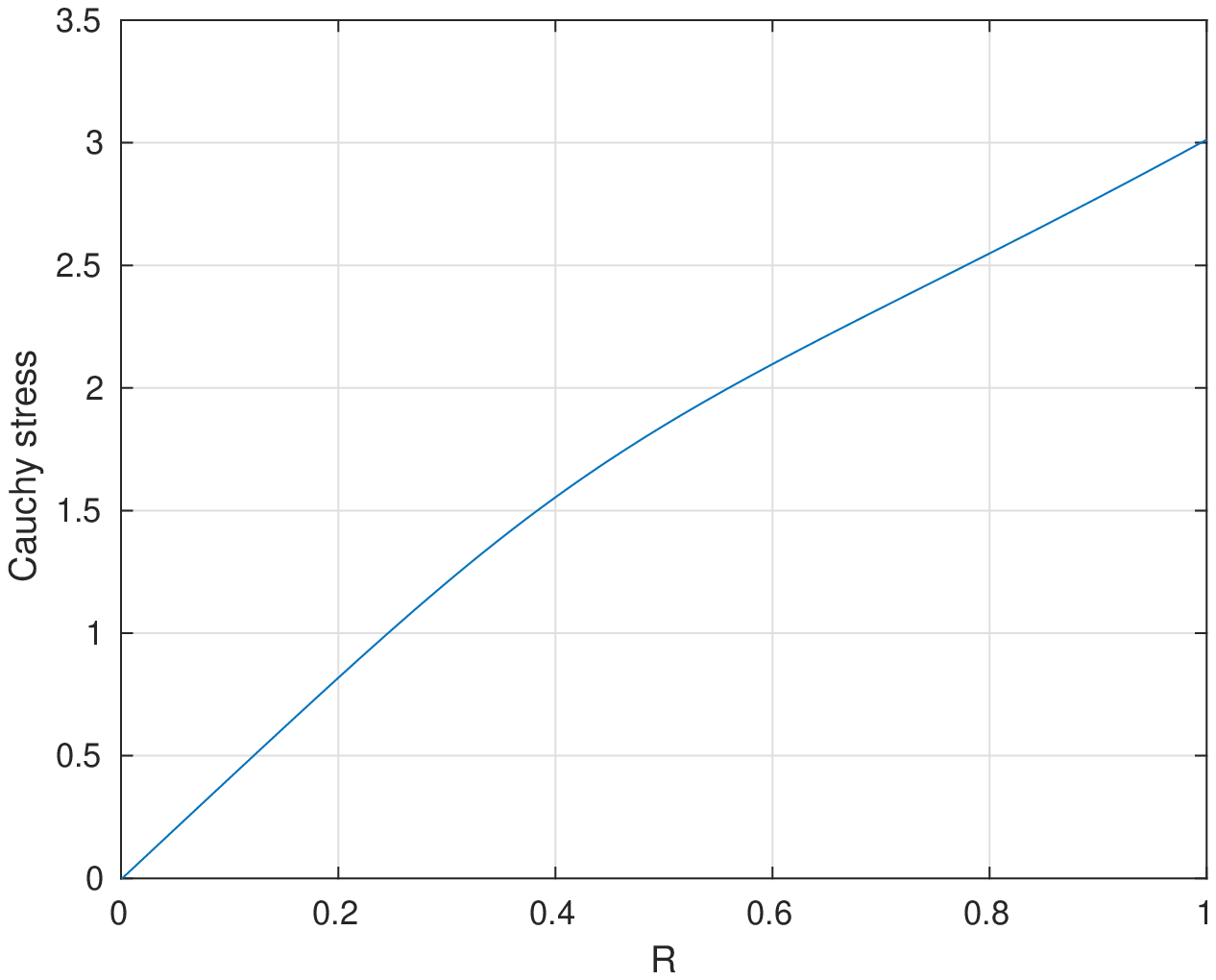}}
\end{center}
\caption{Cauchy stress for $\lambda=1.15$ and $\rho_0=1$.}
\label{fig:8}
\end{figure}

In the our last simulation we take $\lambda=1.15$ and $\rho_0=1$. In Figure 
\ref{fig:6} we show the computed minimizer $r$ which corresponds to a 
cavitating solution. The value of $r(0.001)$ is $0.48346$ with an energy of 
$0.91034$. Figure \ref{fig:7} shows plots of the strains $\deru{r}(R)$ and 
$\frac{r(R)}{R}$, and 
the determinant $\deru{r}(R)(r(R)/R)^2$, and in Figure \ref{fig:8} we show the 
graph of the corresponding Cauchy stress. The boundary layer in $r(R)/R$ 
is now a ``true'' one associated with the computed cavitating solution.
\section{Final comments}
The usual self gravitating problem is that in which the centre of the body 
remains intact ($r(0)=0$) and no condition is explicitly prescribed on the 
outer boundary. This is the problem considered in \cite{JiKoChGo2019} and is a 
special case of one of the problems treated in \cite{CaLe2012}. It is 
straightforward to check that our results hold in this case as well where the 
admissible set is now given by
\[
\mathcal{A}=\big\{r\in
W^{1,1}(0,1)\,|\,r(0)=0,\,\deru{r}(R)>0\mbox{~a.e.
for~}R\in(0,1),\,I_{\mbox{mec}}(r)<\infty\big\}.
\]
In  particular minimizers exist for all densities $\rho_0$ and satisfy the 
Euler--Lagrange equation \eqref{eqn2.13} with natural boundary condition at 
$R=1$ given by $T(1,r(1))=0$ (cf. \eqref{Cstress}). In reference to 
\eqref{se_collapse}. the growth condition \eqref{hdet_growth} with H3 places 
our stored energy function under the ``strong'' compressibility category. 
Thus our result on the existence of minimizers with the centre intact for all 
densities $\rho_0$, is consistent with the results in \cite{JiKoChGo2019} which 
in turn 
suggest that there might be uniqueness of solutions of \eqref{eqn2.13} in 
our problem as well.
\appendix
\section{The gravitational potential}\label{gavpot}
In this section we show how the potential term in the energy functional 
\eqref{eqn2.10} is obtained or follows from the corresponding three dimensional 
potential energy functional.

Let $\B$ be the unit ball with centre at the origin and $\rho_0$ be the mass 
density (per volume) in the reference configuration. We get an expression for 
the potential given in \cite{CaLe2012} by\footnote{The gravitational constant 
$G$ in this expression has been normalized to $\frac{1}{2}$.}
\[
 V(\ts{u})=\frac{1}{2}\int_{\B}\int_{\B}\dfrac{\rho_0(\ts{x})\rho_0(\ts{y})}
 {\norm{\ts{u}(\ts{x})-\ts{u}(\ts{y})}}\,\dif\ts{y}\dif\ts{x},
\]
when $\ts{u}$ is radial (cf. \eqref{eqn2.8}) and $\rho_0$ is radial as well. By 
symmetry, we can set the vertical or $z$ axis in the inner integral to be along 
the direction of $\ts{u}(\ts{x})$. Then by considering a triangle with sides 
$\norm{\ts{u}(\ts{x})}=r(R)$, $\norm{\ts{u}(\ts{y})}=r(U)$, and 
$\norm{\ts{u}(\ts{x})-\ts{u}(\ts{y})}$, we get that
\[ 
\norm{\ts{u}(\ts{x})-\ts{u}(\ts{y})}=\left[r(R)^2+r(U)^2
-2r(R)r(U)\cos\phi\right]^{\frac{1}{2}},
\]
where $R=\norm{\ts{x}}$, $U=\norm{\ts{y}}$, and $\phi$ is the angle between the 
vertical direction (along $\ts{u}(\ts{x})$) and $\ts{u}(\ts{y})$. Using these, 
and writing $\rho_0(\norm{\ts{x}})$ for $\rho_0(\ts{x})$, we get that 
$V(\ts{u})=V(r)$ where
\begin{eqnarray*}
 2V(r)&=&4\pi\int_0^1\rho_0(R)R^2\left[2\pi\int_0^1\int_0^\pi 
\dfrac{\rho_o(U)\sin\phi\,U^2} 
 {\left[r(R)^2+r(U)^2-2r(R)r(U)\cos\phi\right]^{\frac{1}{2}}} 
\,\dif\phi\dif U\right]\,\dif R,\\
&=&4\pi\int_0^1\rho_0(R)R^2\left[2\pi\int_0^1\dfrac{\rho_o(U)U^2} {r(R)r(U)} 
\left(r(R)+r(U)-\abs{r(R)-r(U)}\right)\,\dif U\right]\,\dif R,\\
&=&4\pi\int_0^1\rho_0(R)R^2\left[2\pi\int_0^R\dfrac{\rho_o(U)U^2} {r(R)r(U)} 
2r(U)\,\dif U+2\pi\int_R^1\dfrac{\rho_o(U)U^2} {r(R)r(U)} 
2r(R)\,\dif U\right]\,\dif R,\\
&=&4\pi\int_0^1\rho_0(R)R^2\left[\dfrac{M_R}{r(R)}+4\pi\int_R^1\dfrac{
\rho_o(U)U^2} {r(U)}\,\dif U\right]\,\dif R.
\end{eqnarray*}
Integrating by parts we get that 
\[
 \int_0^1\rho_0(R)R^2\left[4\pi\int_R^1\dfrac{
\rho_o(U)U^2} {r(U)}\,\dif U\right]\,\dif R=
\int_0^1\dfrac{\rho_0(R)M_R}{r(R)}R^2\,\dif R.
\]
Hence 
\[
 V(r)=4\pi\int_0^1\dfrac{\rho_0(R)M_R}{r(R)}R^2\,\dif R.
\]

\pagebreak

\end{document}